# Picone identities for half-linear differential equations of fourth order

**L.M. Cooks[1*),]  Y.A. Stepanyants[2]**

[1] University of New South Wales, Sydney, Australia. E-mail: LevCooks@hotmail.com

[2] Australian Nuclear Science and Technology Organisation, PMB 1, Menai (Sydney), NSW, 2234, Australia. E-mail: Yury.Stepanyants@ansto.gov.au

**Abstract**

Picone-type identities are established for half-linear ODEs of fourth order (one-dimensional $p$-biLaplacian). It is shown that in the linear case they reduce to the known identities for fourth order linear ODEs. Picone-type identity known for two half-linear second-order equations is also generalised to set of $N \geq 3$ equations.

Key words: Picone-type identity, half-linear differential equations, p-Laplacian, Sturmian comparison theorem.

___________________________________

[*)] The author is also known by the name of L.M. Kuks (see [11, 12]).

## 1. Introduction

Almost 100 years ago M. Picone published a paper [1] in which he made a generalisation of the Sturm comparison theorem for the solutions $u$ and $v$ of the pair of ordinary differential equations (ODEs):

$$[p(x)u']' + q(x)u = 0, \tag{1.1}$$

$$[P(x)v']' + Q(x)v = 0. \tag{1.2}$$

Picone discovered that if $u(x)$ and $v(x)$ are differentiable with $p(x)u'$, $P(x)v'$ and $v(x) \neq 0$ on the interval $I \subset \mathbb{R}$, then on this interval the identity holds:

$$\frac{d}{dx}\left\{\frac{u}{v}[vpu' - uPv']\right\} = (p - P)(u')^2 + (Q - q)u^2 + P\left(u' - \frac{u}{v}v'\right)^2. \tag{1.3}$$

Since then this identity bears Picone's name and got various generalisations not only for ordinary but also for elliptic partial differential equations, strongly elliptic systems, equations and systems of even order, time-dependent (parabolic) equations and systems as well as difference equations and systems. The Picone identity turned out to be a useful tool in the development of classical Sturm theorems and study of oscillation and non-oscillation theory of linear differential equations.

Recently there has been an increasing interest in studying the Picone identities for second order half-linear ODEs and $p$-Laplacians (see, e.g., [2–6] and references therein). The main results obtained may be summarised as follows.

Define the function $\varphi(s) = |s|^{\alpha-1}s$, $\alpha > 0$, and consider the differential operators of the form:

$$l[u] \equiv \left[p(x)\varphi(u')\right]' + q(x)\varphi(u), \tag{1.4}$$

$$L[v] \equiv \left[P(x)\varphi(v')\right]' + Q(x)\varphi(v), \tag{1.5}$$

where $p(x)$, $q(x)$, $P(x)$ and $Q(x)$ are continuous functions on a given interval $I \subset \mathbb{R}$. If $u(x)$ and $v(x)$ are differentiable functions with $p(x)\varphi(u')$, $P(x)\varphi(v')$ and $v(x) \neq 0$ on the interval I, then the identity holds in this interval:

$$\frac{d}{dx}\left\{\frac{u}{\varphi(v)}\left[\varphi(v)p\varphi(u') - \varphi(u)P\varphi(v')\right]\right\} = (p-P)(u')^{\alpha+1} + (Q-q)|u|^{\alpha+1} + P\left[|u'|^{\alpha+1} + \alpha\left|u\frac{v'}{v}\right|^{\alpha+1} - (\alpha+1)u'\varphi\left(u\frac{v'}{v}\right)\right] + \frac{u}{\varphi(v)}\left\{\varphi(v)l[u] - \varphi(u)L[v]\right\}. \tag{1.6}$$

In the case $\alpha = 1$, $\varphi(s) = s$ the operators (1.4) and (1.5) are linear; the term in square brackets in the right-hand side of (1.6) reduces to

$$|u'|^{\alpha+1} + \alpha\left|u\frac{v'}{v}\right|^{\alpha+1} - (\alpha+1)u'\varphi\left(u\frac{v'}{v}\right) = (u')^2 + \left(u\frac{v'}{v}\right)^2 - 2u'u\frac{v'}{v} = \left(u' - \frac{u}{v}v'\right)^2, \tag{1.7}$$

and the identity (1.6) coincides with (1.3).

Several kinds of Picone-type identities are known for linear fourth order ODEs and systems of fourth order ODEs [7–10]. Further generalisation can be done for half-linear ODEs of fourth order. One of the objectives of this paper is to establish the Picone-type identities for half-linear

ODEs of fourth order (one-dimensional *p*-biLaplacian). Another objective is to generalise the results known for two half-linear second-order equations to set of $N \geq 3$ equations.

**2. Main results**

*2.1. Picone-type identities for half-linear ODEs of fourth order.*

Consider a pair of equations

$$l[u] \equiv [a(x)\varphi(u'')]'' - [b(x)\varphi(u')]' + c(x)\varphi(u), \tag{2.1}$$

$$L[v] \equiv [A(x)\varphi(v'')]'' - [B(x)\varphi(v')]' + C(x)\varphi(v). \tag{2.2}$$

By the straightforward differentiation the identities following from these equations can be proven.

**Identity 1**: If $v(x)$ does not vanish in the interval I, then the identity holds on this interval:

$$\begin{aligned}
&\frac{d}{dx}\left\{\left[u(a\varphi(u''))' - a\varphi(u'')u'\right] + \left[A\varphi(v'')\left(\frac{u\varphi(u)}{\varphi(v)}\right)' - \frac{u\varphi(u)}{\varphi(v)}(A\varphi(v''))'\right] - \right.\\
&\left.\frac{u}{\varphi(v)}\left[\varphi(v)b\varphi(u') - \varphi(u)B\varphi(v')\right]\right\} = \frac{u}{\varphi(v)}\{\varphi(v)l[u] - \varphi(u)L[v]\} + \\
&A\alpha(\alpha+1)|u|^{\alpha-1}\varphi\left(\frac{v''}{v}\right)\left(\frac{u'v - uv'}{v}\right)^2 - A\left[|u''|^{\alpha+1} + \alpha\left|\frac{u}{v}v''\right|^{\alpha+1} - (\alpha+1)u''\varphi\left(\frac{u}{v}v''\right)\right] - \\
&B\left[|u'|^{\alpha+1} + \alpha\left|\frac{u}{v}v'\right|^{\alpha+1} - (\alpha+1)u'\varphi\left(\frac{u}{v}v'\right)\right] + (A-a)|u''|^{\alpha+1} + (B-b)|u'|^{\alpha+1} + (C-c)|u|^{\alpha+1}.
\end{aligned} \tag{2.3}$$

**Identity 2**: If none of $v(x)$ and $v'(x)$ vanish in the interval I, then the identity holds on this interval:

$$\frac{d}{dx}\left\{\frac{u}{\varphi(v)}\left[\varphi(u)(A\varphi(v''))' - \varphi(v)(a\varphi(u''))'\right] + \frac{u'}{\varphi(v')}\left[\varphi(v')(a\varphi(u''))' - \varphi(u')(A\varphi(v''))'\right] + \frac{u}{\varphi(v)}\left[\varphi(v)b\varphi(u') - \varphi(u)B\varphi(v')\right]\right\} = (a-A)|u''|^{\alpha+1} + (b-B)|u'|^{\alpha+1} + (c-C)|u|^{\alpha+1} +$$

$$A\left[|u''|^{\alpha+1} + \alpha\left|u'\frac{v''}{v'}\right|^{\alpha+1} - (\alpha+1)u''\varphi\left(u'\frac{v''}{v'}\right)\right] + \quad (2.4)$$

$$\left[B|v'|^{\alpha+1} - v'(A\varphi(v''))'\right]\left[\left|\frac{u'}{v'}\right|^{\alpha+1} + \alpha\left|\frac{u}{v}\right|^{\alpha+1} - (\alpha+1)\frac{u'}{v'}\varphi\left(\frac{u}{v}\right)\right] +$$

$$\frac{u}{\varphi(v)}\{\varphi(u)L[v] - \varphi(v)l[u]\}.$$

In the linear case ($\alpha = 1$, $\varphi(s) = s$) the identity (2.3) coincides with the Picone-type identity which has been established by Dunninger [7, 8], and the identity (2.4) coincides with the Picone-type identity established by Kusano & Yoshida [9].

Next two comparison theorems can be proven by using identities (2.3) and (2.4) similarly to how it was done for the second-order half-linear differential equations [3–6].

**Theorem 1**: Suppose $u(x)$ is a nontrivial solution of the equation $l[u] = 0$ in $(a, b)$ such that $u(a) = u(b) = u'(a) = u'(b) = 0$ and $0 \leq A(x) \leq a(x)$, $0 \leq B(x) \leq b(x)$, $C(x) \leq c(x)$ in $[a, b]$. Then, every solution $v(x)$ of the equation $L[v] = 0$ in $(a, b)$ which satisfies $v''(x)/v(x) < 0$ in $[a, b]$ has a zero on $[a, b]$ or is a constant multiple of $u(x)$.

In the linear case this theorem has been proven by Dunninger [7, 8].

**Theorem 2**: Assume that $A(x) \geq 0$ in $I = (a, b)$. If there exists a nontrivial solution $u(x)$ of the equation $l[u] = 0$ in $(a, b)$ such that $u(a) = u(b) = u'(a) = u'(b) = 0$ and $A(x) \leq a(x)$, $B(x) \leq b(x)$, $C(x) \leq c(x)$ in $[a, b]$, then, every solution $v(x)$ of the equation $L[v] = 0$ in $(a, b)$ which satisfying to the inequality $B|v'|^{\alpha+1} - v'(A\varphi(v''))' > 0$ in $(a, b)$ has a zero on $[a, b]$ or is a constant multiple of $u(x)$.

In the linear case this theorem has been proven in [9, 10].

*2.2. Picone-type identities for $N \geq 3$ half-linear second equations.*

Now we generalise the Picone-type identity from two differential equations to $N$ such equations, where $N \geq 3$, in the same manner as it was done by one of the authors in [12].

Consider $N$ differential equations of second order

$$L_k[u_k] \equiv (p_k(x)\varphi(u'_k))' + q_k(x)\varphi(u_k) = 0, \quad k = 1, 2, \ldots, N, \tag{2.5}$$

where $p_k(x)$, $q_k(x)$ are continuous functions on a given interval $I \subset \mathbb{R}$. If $u_k(x)$ are solutions of Eqs. (2.5) differentiable with $p_k(x)\varphi(u'_k)$ and not vanish in the interval I, then the following identity holds on this interval:

$$\frac{d}{dx}\left\{u_{N-1}\varphi(u_{N-1})\left[\sum_{k=0}^{N-1}(-1)^{N-k-1}C_{N-1}^k p_{k+1}(x)\frac{\varphi(u'_{k+1})}{\varphi(u_{k+1})}\right]\right\} =$$
$$\left[\sum_{k=0}^{N-1}(-1)^{N-k-1}C_{N-1}^k p_{k+1}(x)\right]u'_{N-1}\varphi(u'_{N-1}) - \left[\sum_{k=0}^{N-1}(-1)^{N-k-1}C_{N-1}^k q_{k+1}(x)\right]u_{N-1}\varphi(u_{N-1}) - \tag{2.6}$$
$$\sum_{k=0}^{N-1}(-1)^{N-k-1}C_{N-1}^k p_{k+1}(x)\left[|u'|_{N-1}^{\alpha+1} + \alpha\left|u'_{k+1}\frac{u_{N-1}}{u_{k+1}}\right|^{\alpha+1} - (\alpha+1)u'_{N-1}\varphi\left(u_{N-1}\frac{u'_{k+1}}{u_{k+1}}\right)\right].$$

This identity also may be verified by the straightforward differentiation. If $N = 2$, the identity reduces to (1.6). If $N = 3$, then the following comparison theorem can be proven:

**Comparison theorem**: Suppose $u_2(x)$ is a nontrivial solution of the equation $L_2[u_2] = 0$ in $[a, b]$ such that $u_2(a) = u_2(b) = 0$ and $p_k(x) > 0$, $p_2(x) \geq \frac{1}{2}[p_1(x) + p_3(x)]$, $q_2(x) \leq \frac{1}{2}[q_1(x) + q_3(x)]$ in $[a, b]$. Then, every solution $u_1(x)$ or $u_3(x)$ of the equation $L_1[u_1] = 0$ and $L_3[u_3] = 0$ has a zero on $[a, b]$ or is a constant multiple of $u_2(x)$.

Identities (2.3) and (2.4) also can be extended in such way. The results obtained can be also generalised for half-linear PDEs.